\documentclass[12pt, leqno]{amsart}
\usepackage[mathscr]{eucal}
\usepackage{amssymb,amsmath}
\usepackage{amsfonts}
\usepackage{a4}

\newtheorem{theorem}{Theorem}

\newtheorem{lemma}{Lemma}
\newtheorem{corollary}{Corollary}

\newtheorem{definition}{Definition}
\begin{document}
\title{An infinitesimal-birational duality through differential operators}
\author[Tomasz Maszczyk]{Tomasz Maszczyk\dag}

\address{Institute of Mathematics\\
Polish Academy of Sciences\\
Sniadeckich 8\newline 00--956 Warszawa, Poland\\
\newline Institute of Mathematics\\
University of Warsaw\\ Banacha 2\newline 02--097 Warszawa, Poland}
\email{maszczyk@mimuw.edu.pl}
\thanks{\dag The author was partially supported by KBN grants 1P03A 036 26 and 115/E-343/SPB/6.PR UE/DIE 50/2005-2008.}

\begin{abstract}
The structure of filtered algebras of Grothendieck's differential
operators of truncated polynomials in one variable and graded
Poisson algebras of their principal symbols is explicitly
determined. A related infinitesimal-birational duality realized by
a Springer type resolution of singularities and the Fourier
transformation is presented. This algebro-geometrical duality is
quantized in appropriate sense and its quantum origin is
explained.
\end{abstract}

\maketitle

\vspace{1cm}
\begin{flushright}\small{\em It seems that our ship is
arriving at a coast,} \\ \small{\em but the land is still hidden
in the fog.}
\end{flushright}

\vspace{1cm}
\paragraph{\bf 1. Introduction.}
The study of the algebra of differential operators of a
commutative algebra was initiated by Grothendieck \cite{grot} in
1967. The most important question in this theory is how properties
of differential operators reflect the properties of a commutative
algebra. Grothendieck proved that the algebra of differential
operators on a smooth affine variety over a field of
characteristic zero is generated by operators of order at most one
\cite{grot}, in particular is finitely generated. The converse,
that generation by first order filtration implies smoothness, is
now referred to as Nakai's conjecture. It implies the
Zariski-Lipman conjecture \cite{lip}, \cite{zar-lip1},
\cite{zar-lip2}. Nakai's conjecture and its variants has been
verified only for few classes of commutative algebras in
\cite{nak-conj}, \cite{nak-conj2}, \cite{nak-conj3},
\cite{nak-conj4}, \cite{nak-conj5}, \cite{nak-cur}, \cite{nakai2},
\cite{nak-norm}. In general, the determination of the structure of
differential operators of a given non-smooth commutative algebra
is a hard problem and one expects some pathology caused by
singularities. For example, the algebra of differential operators
on the cubic cone $x^{3} + y^{3} +  z^{3} = 0$ is not generated by
operators of bounded order, in particular is not finitely
generated \cite{cub}. Although for reduced curves the algebra of
differential operators is Noetherian \cite{muh},
\cite{smith-staff}, for the non-reduced commutative algebra of
Krull dimension one,  $k[x, y]/(x^{2}, xy)$, its algebra of
differential operators is right but not left Noetherian
\cite{muh}. Such and other problems concerning the structure of
the algebra of differential operators can be found also in
\cite{brown1}, \cite{brown2}, \cite{brown3}, \cite{stan-rei},
\cite{jon}, \cite{levas1}, \cite{levas2}, \cite{muh}, \cite{mus1},
\cite{mus2}, \cite{sai1}, \cite{sai2}, \cite{prime},
\cite{smith-staff}, \cite{mon}, \cite{stan-rei2}.

In this paper we determine the structure of the filtered algebra
${\mathcal D}({\mathcal O}_{n})$ of differential operators of the
polynomial algebra ${\mathcal O}_{n}=k[x]/(x^{n+1})$ of the $n$-th
infinitesimal neighborhood of a point on the affine line and its
associated graded Poisson algebra ${\mathcal P}({\mathcal
O}_{n})$. It turns out that instead of expected pathology we find
beautiful and intriguing phenomena.

Let $(e, h, f)$ be the standard basis of ${\rm sl}_{2}$, subject
to the relations: $[e, f] =  h,\ [h, e] =  2e,\ [h, f]  =  -2f$.
We show that there exists an extension of filtered algebras
$$(e^{n+1}) \rightarrowtail {\mathcal U}({\rm sl}_{2})/(C -
n(n+2))\twoheadrightarrow {\mathcal D}({\mathcal O}_{n}),$$ where
$C=h^{2}+2(ef+fe)$ is the Casimir element and we filter the
enveloping algebra giving $(e, h, f)$ orders $(0, 1, 2)$. The
epimorphism maps $(e, h, f)$ onto operators
$$(x,\ \   2x\frac{{\rm d}}{{\rm d}x} - n,\ \  -x\frac{{\rm d}^{2}}{{\rm d}x^{2}} + n\frac{{\rm d}}{{\rm d}x}).$$
In particular, ${\mathcal D}({\mathcal O}_{n})$ for every $n$ is
generated by operators of order less or equal 2, or more
precisely, by $x$ and a single operator of order 2. This means
that algebras ${\mathcal O}_{n}$ of $n$-th infinitesimal
neighborhoods of a one-dimensional smooth point are close to
smooth algebras.

After passing to the associated graded Poisson algebras we get the
extension
$$\mathfrak m^{n+1} \rightarrowtail {\mathcal P} \twoheadrightarrow {\mathcal P}({\mathcal O}_{n}),$$
where ${\mathcal P}$ is the algebra of polynomial functions on the
nilpotent cone in ${\rm sl}_{2}^{*}$ with its canonical graded
Kirillov-Kostant Poisson structure and $\mathfrak m$ is the
maximal ideal of the vertex. Both extensions are ${\rm
sl}_{2}$-invariant in an appropriate way. In the second extension
the epimorphism describes the closed embedding via the moment map
associated with a hamiltonian action of ${\rm sl}_{2}$ on
Spec$({\mathcal P}({\mathcal O}_{n}))$.

Note that the first extension can be regarded as a quantization of
the second extension of their associated commutative graded
Poisson algebras. Although maximal primitive quotients of the
enveloping algebra ${\mathcal U}({\rm sl}_{2})/(C - n(n + 2))$
depend on $n$, they are Morita equivalent one to each other by the
Beilinson-Bernstein theorem \cite{bei-ber}. This means that the
Morita equivalence class of ${\mathcal U}({\rm sl}_{2})/(C - n(n +
2))$, a quantization of the nilpotent cone in ${\rm sl}_{2}^{*}$,
does not depend on $n$. It is easy to see that forgetting about
filtration one gets ${\mathcal D}({\mathcal O}_{n})\cong
\mbox{End}_{k}({\mathcal O}_{n})$ so for all $n$ algebras
${\mathcal D}({\mathcal O}_{n})$ are Morita equivalent. This means
that for all $n$ the above quantizations of $n$-th infinitesimal
neighborhoods of the vertex of the nilpotent cone are equivalent.

On the other hand, although for $n>0$ surjective homomorphisms of
commutative algebras ${\mathcal O}_{n+1}\rightarrow{\mathcal
O}_{n}$ do not induce homomorphisms of algebras of differential
operators ${\mathcal D}({\mathcal O}_{n+1})\rightarrow{\mathcal
D}({\mathcal O}_{n})$ (there is no any homomorphism of algebras of
matrices ${\mathcal M}_{n+1}\rightarrow {\mathcal M}_{n}$), we
have surjective homomorphisms of their graded Poisson algebras
${\mathcal P}({\mathcal O}_{n+1})\rightarrow{\mathcal P}({\mathcal
O}_{n})$ which is a kind of surprise, which cannot be explained by
classical geometry.

Next, after passing to the inverse limits, we compare the
completion of the algebra of principal symbols on the affine line
with the inverse limit of the above system of algebras of
principal symbols on infinitesimal neighborhoods of a point in the
affine line. We show that there exists a unique grading preserving
homomorphism of graded Poisson algebras over the completed local
ring of a point in the line, between these two limits.
Geometrically, this is a completion of a resolution of
singularities in the category of conical Poisson varieties
$$\begin{array}{ccccc}
{\rm T}^{*}\mathbb{A}^{1}     &           &  \stackrel{\exists!}{\longrightarrow} &           & Y \\
                              & \searrow  &                                       & \swarrow  &   \\
                              &           &            \mathbb{A}^{1}             &           &
     \end{array}.$$
Finally, we show that the above canonical resolution of
singularity together with the Fourier transformation establish a
\textbf{duality} between the above $\rm{SL}_{2}$-symmetries of
principal symbols of the system of infinitesimal neighborhoods of
finite order of a point in the line $\mathbb{A}^{1}$
(\textbf{infinitesimal} picture)  and automorphisms of dual
principal symbols induced by birational automorphisms of the dual
line $\hat{\mathbb{A}}^{1}$ (\textbf{birational } picture).

\vspace{3mm}

{\em {\rm \textbf{``Magic Diamond Theorem"}}. The following
diagram of canonical morphisms in the category of conical Poisson
varieties
$$\begin{array}{ccc}
                                              &  {\rm T}^{*}\mathbb{A}^{1}  &              \\
\stackrel{\stackrel{\mbox{\scriptsize{\rm resol.\ sing.}}}{}}{}\swarrow  &                             & \searrow\stackrel{\stackrel{\mbox{\scriptsize{\rm cl. Fourier\ tr.}}}{}}{}     \\
   Y          &                             &   {\rm T}^{*}\hat{\mathbb{A}}^{1}        \\
                                     \mbox{\scriptsize{\rm inf. moment\ m.}}\searrow &                             & \swarrow\mbox{\scriptsize{\rm bir. moment\ m.}}     \\
                                              &   {\rm sl}_{2}^{*}          &
                                   \end{array}$$
is commutative.}

\vspace{3mm}

This diagram can be quantized in appropriate sense, which recovers quantum nature of classical algebraic geometry as follows

\vspace{3mm}

{\em {\rm \textbf{``Quantized Magic Diamond Theorem"}}.  There exists the
following commutative diagram  in the category of
almost commutative algebras}
$$\begin{array}{ccc}
                                              &  \mathcal{D}(\mathbb{A}^{1})  &              \\
\stackrel{\stackrel{\mbox{\scriptsize{\rm q. resol.\ sing.}}}{}}{}\nearrow  &                             & \nwarrow\stackrel{\stackrel{\mbox{\scriptsize{\rm q. Fourier\ tr.}}}{}}{}     \\
  \mathcal{U}({\rm sl}_{2})/(C-n(n+2))   &                             &    \mathcal{D}(\hat{\mathbb{A}}^{1})        \\
                                     \mbox{\scriptsize{\rm q. inf. moment\ m.}}\nwarrow &                             & \nearrow\mbox{\scriptsize{\rm q. bir. moment\ m.}}     \\
                                              &   \mathcal{U}({\rm sl}_{2})          &
                                   \end{array}$$
{\em which is a quantization of the Magic Diamond.}

\vspace{3mm}

It should be understood as a purely mathematical analog of the
physical de Broglie particle-wave duality \cite{Fey}, where an
infinitesimal neighbourhood of a closed point on the line and the
general point of the dual line are mathematical counterparts of a
spatially localized particle and its Fourier decomposition into
waves of frequencies spread in continuous spectrum, respectively.
The de Broglie particle-wave duality could not be explained in
terms of classical physics and was the main motivation for the new
quantum physics. Our infinitesimal-birational duality cannot be
explained in terms of classical geometry (what we emphasize in the
name of the above theorem). In both dualities the non-classical
part is the Fourier transformation realizing both dualities in a
similar way. This suggests that classical geometry should be
``quantized" in an appropriate sense in order to explain the
perfect harmony of the following aspects of our example:
\begin{itemize}
\item $\mathcal{D}$-modules on non-reduced schemes.
\item Representation theory and moment maps.
\item Resolution of singularities.
\item Fourier transformation.
\end{itemize}
It should be noted that our result is not an isolated ``quantum"
phenomenon relating algebraic geometry, representation theory,
moment maps and the Fourier transformation. One attempt in this
direction is the general program of deformation quantization of
nilpotent coadjoint orbits of complex semisimple Lie groups
\cite{AsBry1}, \cite{AsBry2}, \cite{BryKos}, \cite{Dit},
\cite{FioLle}, and deformation quantization of moment maps
\cite{Ham1}, \cite{Ham2}, \cite{Xu}. Another program concentrates
around Vogan's conjecture relating  the Fourier transform of a
nilpotent orbit of the coadjoint representation of a complex Lie
group and multiplicities of the ring of polynomial functions on
that orbit \cite{Kin1}, \cite{Kin2}, \cite{Kin3}, \cite{Ver2}.
Moreover, the Fourier transformation of sheaves supported on the
nilpotent cone in the co-adjoint representation of a reductive
algebraic group play an important role in the theory of Springer
representations of the Weyl group on intersection cohomology
\cite{BorMac}, \cite{Bry}, \cite{Gin}, \cite{Hot}, \cite{KazLus},
\cite{Spr1}, \cite{Spr2}, \cite{Spr3}.

\vspace{3mm}
\paragraph{\textit{Acknowledgement}.} I would like to express
my gratitude to Mariusz Wodzicki for numerous helpful discussions
and encouragement.

\vspace{3mm}
\paragraph{\bf 2. Basic definitions.}

Let $k$ be a field of characteristic zero.
\begin{definition}
A filtered associative algebra ${\mathcal D}=\bigcup_{p\geq
0}{\mathcal D}^{p}$, ${\mathcal D}^{p}\subset{\mathcal D}^{p+1}$,
is called \textbf{almost commutative} if $[{\mathcal D}^{p},
{\mathcal D}^{q}]\subset {\mathcal D}^{p+q-1}$.
\end{definition}
\begin{definition}
Given a commutative algebra ${\mathcal O}$ one defines the
filtered \textbf{algebra}  ${\mathcal D}({\mathcal
O})=\bigcup_{p\geq 0}{\mathcal D}^{p}({\mathcal O})$ \textbf{of
differential operators of} ${\mathcal O}$, where
$${\mathcal D}^{p}({\mathcal O}):=\{\delta\in{\rm End}_{k}({\mathcal O})\mid \forall_{f_{0},\ldots,f_{p}\in{\mathcal O}}\ \
 [f_{0},\ldots,[f_{p},\delta]\ldots]=0\}.$$
This filtration is called \textbf{order filtration}.
\end{definition}

\vspace{3mm}\paragraph{One shows that the filtered associative
algebra ${\mathcal D}({\mathcal O})$ is  almost commutative.}
\begin{definition}
A  graded commutative algebra ${\mathcal P}=\bigoplus_{p\geq
0}{\mathcal P}^{p}$ is called \textbf{graded Poisson algebra} if
there is given a Lie algebra structure on ${\mathcal P}$
$$\{ -, -\} : {\mathcal P}\otimes{\mathcal P}\rightarrow{\mathcal P},$$
such that
$$\{{\mathcal P}^{p}, {\mathcal P}^{q}\} \subset{\mathcal
P}^{p+q-1}$$
and for all $f, g, h\in {\mathcal P}$
$$\{ f, gh\}=\{ f, g\}h+g\{ f, h\}.$$
\end{definition}

\vspace{3mm} \paragraph{\textbf{Example.}}The associated graded
commutative algebra ${\rm Gr}({\mathcal D}):=\bigoplus_{p\geq
0}{\mathcal D}^{p}/{\mathcal D}^{p-1}$ of an almost commutative
filtered algebra ${\mathcal D}$ is a graded Poisson algebra, where
for all $\delta\in{\mathcal D}^{p}$, $\epsilon\in{\mathcal D}^{q}$
$$\{\delta+{\mathcal D}^{p-1}, \epsilon+ {\mathcal D}^{q-1}\}:= [\delta, \epsilon]+{\mathcal D}^{p+q-2}.$$

\begin{definition}
Homogeneous elements of a  graded Poisson algebra ${\mathcal
P}({\mathcal O}):={\rm Gr}({\mathcal D}({\mathcal O}))$ are called
\textbf{principal symbols of differential operators of} ${\mathcal
O}$.
\end{definition}
\vspace{3mm}
For $X={\rm Spec}({\mathcal O})$ we will write ${\mathcal
D}(X):={\mathcal D}({\mathcal O})$,  ${\mathcal
P}(X):={\mathcal P}({\mathcal O})$.

\vspace{3mm}
\paragraph{\textbf{Example.}} Let $\mathfrak{g}=\bigcup_{p\geq
0 }\mathfrak{g}^{p}$, $\mathfrak{g}^{p}\subset\mathfrak{g}^{p+1}$,
be a filtered Lie algebra, i.e.
$$[\mathfrak{g}^{p},\mathfrak{g}^{q}]\subset\mathfrak{g}^{p+q-1}.$$
Let us give its enveloping algebra ${\mathcal U}(\mathfrak{g})$
the induced filtration. Then ${\mathcal U}(\mathfrak{g})$ becomes
an almost commutative algebra whose associated graded Poisson
algebra $\rm{Gr}({\mathcal U}(\mathfrak{g}))$ is the symmetric
algebra ${\rm S}(\rm{Gr}(\mathfrak{g}))$ of the associated graded
vector space $\rm{Gr}(\mathfrak{g})=\bigoplus_{p\geq
0}\mathfrak{g}^{p}/\mathfrak{g}^{p-1}$, where for generators $X\in
\mathfrak{g}^{p},Y\in \mathfrak{g}^{q}$
$$\{ X + \mathfrak{g}^{p-1}, Y+\mathfrak{g}^{q-1}\}:=[X,Y]+\mathfrak{g}^{p+q-2}.$$
Any ${\rm ad}$-invariant homogeneous ideal in ${\rm
S}(\rm{Gr}(\mathfrak{g}))$ defines a graded Poisson structure on
the factor-algebra. One calls it \textbf{graded Kirillov-Kostant
Poisson structure}.

\begin{definition} Let ${\mathcal P}=\bigoplus_{p\geq
0}{\mathcal P}^{p}$ be a graded Poisson algebra. A
\textbf{quantization} of ${\mathcal P}$ is a pair consisting of
the Morita equivalence class of an almost commutative algebra
${\mathcal D}=\bigcup_{p\geq 0}{\mathcal D}^{p}$ together with an
isomorphism of graded Poisson algebras ${\rm Gr}({\mathcal
D})\rightarrow{\mathcal P}$.
\end{definition}
Note that any almost commutative algebra is a quantization of its
associated graded  Poisson algebra.

\begin{definition}
A morphism of conical Poisson schemes
$${\rm Spec}({\mathcal P})\rightarrow{\rm Spec}({\rm S}(\rm{Gr}(\mathfrak{g})))$$
determined by a homomorphism of graded Poisson algebras
$${\rm S}(\rm{Gr}(\mathfrak{g}))\rightarrow {\mathcal P}$$
is called \textbf{moment map}.
\end{definition}

\begin{definition}
For a given morphism of conical Poisson schemes ${\rm Spec}({\mathcal P}_{1})\rightarrow{\rm Spec}({\mathcal P}_{2})$ its \textbf{quantization}
is a pair consisting of the Morita equivalence class of a
filtration preserving homomorphism of almost commutative algebras
$$ {\mathcal D}_{2}\rightarrow  {\mathcal D}_{1}$$
and an isomorphism of its associated grading preserving
homomorphism of graded Poisson algebras with a given homomorphism of graded Poisson algebras
${\mathcal P}_{2}\rightarrow {\mathcal P}_{1}$. The definition of a quantization of a given diagram in the category of  conical Poisson schemes is obvious.
\end{definition}
Note that the set of $k$-points of ${\rm Spec}({\rm
S}(\rm{Gr}(\mathfrak{g})))$ can be canonically identified with the
dual space $\rm{Gr}(\mathfrak{g})^{*}$. Since all our structures
are defined over $k$ we will abuse the language and will not
distinguish ${\rm Spec}({\rm S}(\rm{Gr}(\mathfrak{g})))$ from
$\rm{Gr}(\mathfrak{g})^{*}$. In the case of a filtered  Lie
algebra coming from a graded Lie algebra
$\mathfrak{g}=\bigoplus_{i\geq 0 }\mathfrak{g}_{i}$,
$[\mathfrak{g}_{i},\mathfrak{g}_{j}]\subset\mathfrak{g}_{i+j-1}$,
i.e. when the filtration has the form
$\mathfrak{g}^{p}=\bigoplus_{0\leq i\leq p }\mathfrak{g}_{i}$, we
will not distinguish also the graded vector space
$\rm{Gr}(\mathfrak{g})$ from the graded Lie algebra
$\mathfrak{g}$.

If one gives a commutative algebra ${\mathcal O}$ the trivial
grading and the trivial Poisson structure then one has the
canonical injective homomorphism of graded Poisson algebras
${\mathcal O}\rightarrowtail {\mathcal P}({\mathcal O})$ onto the
subalgebra of symbols of degree zero, and the canonical surjective
homomorphism of graded ${\mathcal O}$-algebras ${\mathcal
P}({\mathcal O})\twoheadrightarrow{\mathcal O}$ annihilating all
symbols of positive degree. If ${\mathcal O}={\mathcal O}(V)$ is
the algebra of polynomial functions on a smooth affine variety $V$
then ${\rm Spec}({\mathcal P}({\mathcal O}))={\rm T}^{*}V$, the
cotangent bundle of $V$. Then the above homomorphisms describe the
canonical projection ${\rm T}^{*}V\rightarrow V$ and the embedding
onto the zero section $V\rightarrow{\rm T}^{*}V$, respectively.

Let $E$ be a finite dimensional graded vector space and
$\hat{E}=E^{*}[1]$ be its dual with the dual grading enlarged by
one. Then
$${\rm T}^{*}{\rm Spec}({\rm Sym}(E))={\rm Spec}({\rm Sym}(E)\otimes {\rm Sym}(\hat{E})),$$
is a conical Poisson variety with the following Poisson structure
$$\{ e_{1}\otimes 1, e_{2}\otimes 1\} =0,\ \ \{ 1\otimes \hat{e}_{1}, 1\otimes
\hat{e}_{2}\}=0,\ \ \{ 1\otimes \hat{e}, e\otimes
1\}=\hat{e}(e).$$ Using the canonical linear isomorphism of graded
vector spaces $E\rightarrow \hat{\hat{E}}$, we obtain the
canonical identification
$${\rm T}^{*}{\rm Spec}({\rm Sym}(\hat{E}))\stackrel{\cong}{\rightarrow}{\rm Spec}({\rm Sym}(\hat{E})\otimes {\rm Sym}(E)).$$
\begin{definition}
We define an isomorphism of conical Poisson varieties
$${\rm T}^{*}{\rm Spec}({\rm Sym}(E))\rightarrow {\rm T}^{*}{\rm Spec}({\rm Sym}(\hat{E}))$$
by means of polynomial functions as follows
$${\rm Sym}(E)\otimes {\rm Sym}(\hat{E})\leftarrow{\rm Sym}(\hat{E})\otimes {\rm Sym}(E),$$
$$\hat{e}\otimes 1\mapsto 1\otimes \hat{e},\ \ 1\otimes e\mapsto - e\otimes 1.$$
We call this isomorphism \textbf{classical Fourier transformation}.
\end{definition}
Note that our classical Fourier transformation on the level of cotangent
bundles ${\rm T}^{*}(\mathbb{A}^{n})\rightarrow {\rm T}^{*}(\hat{\mathbb{A}}^{n})$ can be regarded as a result of descending of the following  \textbf{quantum Fourier transformation}
\begin{align*}
\mathcal{D}(\hat{\mathbb{A}}^{n}) \rightarrow  \mathcal{D}(\mathbb{A}^{n})
\end{align*}
\begin{align*}
\hat{x}_{i} & \mapsto \frac{\partial}{\partial x_{i}},\ \ \
\frac{\partial}{\partial \hat{x}_{i}} \mapsto - x_{i}
\end{align*}
to the level of principal symbols, what justifies our terminology.

According to Definition 7 the quantum Fourier transformation is a quantization of the classical Fourier transformation, provided $\mathcal{D}(\hat{\mathbb{A}}^{n})$ is filtered not by the order of differential operators but by the degree of their polynomial coefficients. Whenever $\mathcal{D}(\hat{\mathbb{A}}^{n})$ appears in this article it is equipped with this filtration making it an almost commutative algebra.

For the affine line $\mathbb{A}^{1}=\mathbb{P}^{1}\setminus \{
\infty\}$ we can define the dual line as
$\hat{\mathbb{A}}^{1}=\hat{\mathbb{P}}^{1}\setminus \{
\infty\}=\mathbb{P}^{1}\setminus \{ 0\}$, by self-duality of the
projective line, and we regard ${\rm T}^{*}\hat{\mathbb{A}}^{1}$
as an open subscheme in ${\rm T}^{*}\hat{\mathbb{P}}^{1}$. The
affine algebraic group ${\rm SL}_{2}$ acts on
$\hat{\mathbb{P}}^{1}$ and this action lifts canonically to the
conical structure preserving Poisson action on the conical Poisson
variety  ${\rm T}^{*}\hat{\mathbb{P}}^{1}$. We can regard the
restriction of this action to ${\rm T}^{*}\hat{\mathbb{A}}^{1}$ as
the induced rational action of the pseudogroup of birational
automorphisms of $\hat{\mathbb{A}}^{1}$. The Lie algebra ${\rm
sl}_{2}$ of vector fields on $\hat{\mathbb{P}}^{1}$ lifts
canonically to ${\rm T}^{*}\hat{\mathbb{P}}^{1}$ and defines a
moment map ${\rm T}^{*}\hat{\mathbb{P}}^{1}\rightarrow {\rm
sl}_{2}^{*}$. The restriction of this moment map to ${\rm
T}^{*}\hat{\mathbb{A}}^{1}$ can be understood as the moment map
$${\rm
T}^{*}\hat{\mathbb{A}}^{1}\rightarrow {\rm sl}_{2}^{*}$$ defined
by the induced rational action of the pseudogroup of birational
automorphisms of $\hat{\mathbb{A}}^{1}$.

\vspace{3mm}\paragraph{\bf 3. Differential operators of truncated
polynomials.} We will study the polynomial algebra ${\mathcal O}
_{n} := k[x]/(x^{n+1})$ of the $n$-th infinitesimal neighborhood
of a point on the affine line,  the algebra ${\mathcal
D}({\mathcal O} _{n})$ of differential operators on this
neighborhood and the corresponding associated graded Poisson
algebra ${\mathcal P}({\mathcal O}_{n})$ spanned by their
principal symbols. We identify functions from ${\mathcal O}_{n}$
with operators of zero order from ${\mathcal D}^{0}({\mathcal O
}_{n})$ and their principal symbols from ${\mathcal
P}^{0}({\mathcal O }_{n})$, for instance we identify $x^{l}\in
{\mathcal O}_{n} $ with the operator $\delta_{0}^{l}\in {\mathcal
D}^{0}({\mathcal O }_{n})$ of multiplication by $x^{l}$ and with
its principal symbol $y_{0}^{l}\in {\mathcal P}^{0}({\mathcal O
}_{n})$. For any $\delta\in {\mathcal D}({\mathcal O }_{n})$ and
$f\in {\mathcal O}_{n}$ we use the notation
$$\mbox{ad}_{f(x)}(\delta ) := [f(\delta_{0}), \delta ].$$

\begin{lemma}
$${\rm ad}^{p+1}_{x}(\delta  ) = 0 \Leftrightarrow \forall_{f_{0}, ..., f_{p}\in
{\mathcal O}_{n}}\ \ {\rm ad}_{f_{0}} ... {\rm ad}_{f_{p}}(\delta
) = 0.$$
\end{lemma}

\textit{Proof.} One can assume  $f_{0} = x^{l_{0}}, ..., f_{p} =
x^{l_{p}}$, for $l_{0}, ..., l_{p} \geq 0$. \break Since
$\mbox{ad}_{x^{l}}(\delta  ) =
\sum_{m=1}^{l}\delta_{0}^{m-1}\mbox{ad}_{x}(\delta
)\delta_{0}^{l-m}$ and $\mbox{ad}_{x}$ commutes with left and
right multiplication by $\delta_{0}$ we get

$$\mbox{ad}_{x^{l_{0}}}...\mbox{ad}_{x^{l_{p}}}(\delta  ) = $$
$$= \sum_{v_{0}=1}^{l_{0}}\delta_{0}^{v_{0}-1}\mbox{ad}_{x}(\sum_{v_{1}=1}^{l_{1}}\delta_{0}^{v_{1}-1}\mbox{ad}_{x}(...\delta_{0}^{v_{p}-1}\mbox{ad}_{x}(\delta
)\delta_{0}^{l_{p}-v_{p}})...)\delta_{0}^{l_{0}-v_{0}} = $$

$$= \sum_{v_{0}=1}^{l_{0}}...\sum_{v_{p}=1}^{l_{p}}\delta_{0}^{(v_{0}-1)+...+(v_{p}-1)}\mbox{ad}_{x}^{p+1}(\delta)\delta_{0}^{(l_{0}-v_{0})+...+(l_{p}-v_{p})}.$$
which implies the lemma. $\Box $

\begin{lemma}
$$x^{n+1} = 0 \Rightarrow {\rm ad}_{x}^{2n+1} = 0.$$
\end{lemma}

\textit{Proof.} Substituting $l=2n+1$ into the identity
$$\mbox{ad}_{x}^{l}(\delta ) = \sum_{m=0}^{p}(-1)^{m+1}\left( \begin{array}{c} l \\ m
\end{array}\right) \delta_{0}^{l-m}\delta \delta_{0}^{m}$$
one can see that at least one of $m, l-m$ is greater than $n$.
$\Box $

\begin{corollary}
$${\mathcal D}({\mathcal O}_{n}) = {\rm End}_{k}({\mathcal O}_{n}).$$
\end{corollary}

\begin{corollary} The order filtration on differential operators has the
form
$${\mathcal D}^{p}({\mathcal O
}_{n}) = \bigl\{ \delta\in {\rm End}_{k}({\mathcal O}_{n}) \mid
{\rm ad}^{p+1}_{x}(\delta ) = 0 \bigr\} . $$
\end{corollary}

\begin{lemma}
We have the canonical linear embedding ${\mathcal P}^{p}({\mathcal
O }_{n}) \rightarrowtail {\mathcal O}_{n}$ induced by the  map
${\rm ad}^{p}_{x}: {\mathcal D}^{p}({\mathcal O }_{n})\rightarrow
{\mathcal O}_{n} $. It is multiplicative in the following sense:
for all $\delta \in {\mathcal D}^{p}({\mathcal O }_{n})$,
$\epsilon\in {\mathcal D}^{q}({\mathcal O }_{n})$
\begin{align}
\frac{1}{(p+q)!}{\rm ad}^{p+q}_{x}(\delta\epsilon) =
\frac{1}{p!}{\rm ad}^{p}_{x}(\delta)\frac{1}{q!}{\rm
ad}^{q}_{x}(\epsilon).
\end{align}
\end{lemma}

\textit{Proof.} Since ${\mathcal D}^{p}({\mathcal O }_{n}) = \ker
\mbox{ad}^{p+1}_{x}$ then $\mbox{ad}^{p}_{x}$ maps ${\mathcal
D}^{p}({\mathcal O }_{n})$ into ${\mathcal O}_{n} = \ker
\mbox{ad}_{x}$ with kernel ${\mathcal D}^{p-1}({\mathcal O
}_{n})\subset {\mathcal D}^{p}({\mathcal O }_{n})$, hence embeds
${\mathcal D}^{p}({\mathcal O }_{n}) = {\mathcal D}^{p}({\mathcal
O }_{n})/{\mathcal D}^{p-1}({\mathcal O }_{n})$ into ${\mathcal
O}_{n}$. Using the fact that $\mbox{ad}^{l}_{x}(\delta ) = 0$ and
$\mbox{ad}^{m}_{x}(\delta ') = 0$ for $l>p$ and $m>q$ in the
identity
$$\frac{1}{(p+q)!}\mbox{ad}^{p+q}_{x}(\delta\epsilon) =\sum_{l+m=p+q}
\frac{1}{l!}\mbox{ad}^{l}_{x}(\delta)\frac{1}{m!}\mbox{ad}^{m}_{x}(\epsilon)$$
we get on the right hand side only one summand with $l=p$ and
$m=q$. $\Box $

Note that ${\rm ad}^{p}_{x}: {\mathcal D}^{p}({\mathcal O
}_{n})\rightarrow {\mathcal O}_{n} $ is ${\mathcal O}_{n}$-linear,
hence its image is an ideal, necessarily generated by $x^{v_{p}}$
for some $v_{p}\leq n$.

\begin{definition} We define operators $\delta_{p}\in {\mathcal D}^{p}({\mathcal O
}_{n})$ for $p = 1, ..., 2n$ such that
\begin{align}
\frac{1}{p!}{\rm ad}^{p}_{x}(\delta_{p}) = \delta_{0}^{v_{p}}.
\end{align}
\end{definition}

Note that operators $\delta_{p}$ are determined up to ${\mathcal
D}^{p-1}({\mathcal O }_{n})={\rm ker}({\rm ad}^{p}_{x})$. However,
by Lemma 3, their principal symbols $y_{p}\in {\mathcal
P}^{p}({\mathcal O }_{n})$ are defined uniquely.

\begin{lemma} The system $(y_{0}^{l}y_{p})_{l=0}^{n-v_{p}}$ is a basis of
${\mathcal P}^{p}({\mathcal O }_{n})$ for $i = 1, ..., 2n$.
\end{lemma}

\textit{Proof.} Since on the right hand side of the equality
$$\frac{1}{p!}\mbox{ad}^{p}_{x}(\delta_{0}^{l}\delta_{p}) = \delta_{0}^{l+v_{p}}$$
we get linearly independent powers of $\delta_{0}$ then the system
of corresponding principal symbols is linearly independent in
${\mathcal P}^{p}({\mathcal O }_{n})$. On the other hand
 for every $\delta \in {\mathcal D}^{p}({\mathcal O
}_{n})$ we have the following decomposition

$$\mbox{ad}^{p}_{x}(\delta ) = \sum_{l=0}^{n-v_{p}}c_{l}\delta_{0}^{l+v_{p}} =
\sum_{l=0}^{n-v_{p}}c_{l}\delta_{0}^{l} \delta_{0}^{v_{p}} =$$
$$= \sum_{l=0}^{n-v_{p}}c_{l}\delta_{0}^{l}\frac{1}{p!}\mbox{ad}^{p}_{x}(\delta_{p}) =
\mbox{ad}^{p}_{x}(\sum_{l=0}^{n-v_{p}}\frac{c_{l}}{p!}\delta_{0}^{l}\delta_{p})$$
which means that
$$\delta  \equiv \sum_{l=0}^{n-v_{p}}\frac{c_{l}}{p!}\delta_{0}^{l}\delta_{p}\ \
\  {\rm mod}\  {\mathcal D}^{p-1}({\mathcal O }_{n}).$$ Therefore
the above system of principal symbols generates ${\mathcal
P}^{p}({\mathcal O }_{n})$. $\Box $

\begin{corollary}
$$\dim {\mathcal P}^{p}({\mathcal O
}_{n}) = n - v_{p}+1.$$
\end{corollary}

\begin{corollary}
$$v_{1}+...+v_{2n} = n(n+1).$$
\end{corollary}

\textit{Proof.} We have

$$(n+1)^{2} = \dim \mbox{End}_{k}({\mathcal O}_{n}) = \dim {\mathcal D}({\mathcal O}_{n}) = $$
$$= \dim {\mathcal O}_{n} + \sum_{p=1}^{2n}\dim {\mathcal P}^{p}({\mathcal O
}_{n}) = (n+1) + 2n(n+1) - (v_{1} +...+ v_{2n}). \Box $$
 Our next task is to determine $v_{p}$'s for $i = 1, ..., 2n$.

\begin{lemma} $v_{p} > 0$.
\end{lemma}

\textit{Proof.} In the ordered basis $(1, x, ..., x^{n})$ of
${\mathcal O}_{n}$ the multiplication by $x$ has the form of the
sub-diagonal Jordan block. Therefore we have
$${\rm tr}(\delta_{0}^{v_{p}}) = \frac{1}{p!}{\rm tr}(\mbox{ad}^{p}_{x}(\delta_{p})) = 0.$$
But it is possible only if  $ v_{p} > 0$.  $\Box $

\begin{lemma} If $n>0$ then $ v_{1} =  v_{2} = 1$.
\end{lemma}

\textit{Proof.} We know already that $ v_{1},  v_{2} \geq 1$. We
show that $ v_{1},  v_{2} \leq 1$. Using the Jordan form we can
easily find a solution $\delta_{2}\in \mbox{End}_{k}({\mathcal
O}_{n})$ to the equation

\begin{align}\frac{1}{2}\mbox{ad}^{2}_{x}(\delta_{2}) = \delta_{0}.\end{align}
Therefore $v_{2} = 1$.  If we take
$\delta_{1}:=\frac{1}{2}\mbox{ad}_{x}(\delta_{2})$ then
\begin{align}\mbox{ad}_{x}(\delta_{1}) =  \frac{1}{2}\mbox{ad}^{2}_{x}(\delta_{2}) = \delta_{0}.\end{align}
Therefore $v_{1} = 1$.   $\Box $

One can easily check that one can choose the operator $\delta_{2}$
of the form $ x \frac{{\rm d}^{2}}{{\rm d}x^{2}} - n\frac{{\rm
d}}{{\rm d}x}$. Then the operators
\begin{align}
\delta_{0} := x,\ \  \delta_{1} := -x\frac{{\rm d}}{{\rm d}x} +
\frac{n}{2} ,\ \  \delta_{2} := x \frac{{\rm d}^{2}}{{\rm d}x^{2}}
- n\frac{{\rm d}}{{\rm d}x}\end{align}
 satisfy the relations
\begin{align}
[\delta_{0}, \delta_{1}]  =  \delta_{0},  \ \ [\delta_{0},
\delta_{2}] = 2\delta_{1}, \ \ [\delta_{1}, \delta_{2}]  =
\delta_{2},
\end{align}
and
\begin{align}\delta_{1}^{2} -\frac{1}{2}(\delta_{0}\delta_{2} + \delta_{2} \delta_{0}) =
\frac{n}{2}(\frac{n}{2} + 1),\end{align} which means that we have
a homomorphism of filtered algebras from a maximal primitive
quotient of the enveloping algebra of ${\rm sl}_{2}$ \cite{env}
into ${\mathcal D}({\mathcal O}_{n})$
$${\mathcal U}({\rm sl}_{2})/(C - n(n+2))\rightarrow {\mathcal D}({\mathcal O}_{n}),$$
\begin{align}
e\mapsto \delta_{0},\ \ h\mapsto -2\delta_{1},\ \ f\mapsto
-\delta_{2},
\end{align}
where ${\rm sl}_{2}$ is spanned by $(e, h, f)$ subject to the
relations
\begin{align*}
[e, f]  =  h, \ \ [h, e]  =  2e,  \ \ [h, f]  =  -2f,
\end{align*}
and $C$ denotes the Casimir element $h^{2}+2(ef+fe)$.

\begin{definition} We define on the maximal primitive quotient
a filtration coming from the filtration on the Lie algebra ${\rm
sl}_{2}$: $(e)\subset (e, h)\subset (e, h, f)={\rm sl}_{2}$.  We
call the above choice {\rm (5)} of $\delta_{0}, \delta_{1},
\delta_{2}$ and the resulting homomorphism {\rm (8)} of filtered
algebras {\textbf{\textit{distinguished}}}.
\end{definition}

\vspace{3mm} Note that this filtration differs from the standard
filtration of the enveloping algebra of a Lie algebra, under which
the Lie algebra entirely lies in the first piece of the
filtration.

\begin{lemma} If $n>0$ then $ v_{2l-1} =  v_{2l} = l$ for $l = 1, ..., n$.
\end{lemma}

\textit{Proof.} Note first that if $x^{l}\in {\rm
ad}_{x}^{p}({\mathcal D}({\mathcal O}_{n})^{p})\subset {\mathcal
O}_{n} $ then $v_{p}\leq l$. Using the multiplicative law from
Lemma 3 and the equalities $ v_{1} = v_{2} = 1$ we get
\begin{align}
\frac{1}{(2l)!}\mbox{ad}^{2l}_{x}(\delta_{2}^{l}) =
(\frac{1}{2}\mbox{ad}^{2}_{x}(\delta_{2}))^{l} =
\delta_{0}^{l},\end{align}
\begin{align}\frac{1}{(2l-1)!}\mbox{ad}^{2l-1}_{x}(\delta_{1}\delta_{2}^{l-1}) =
\mbox{ad}_{x}(\delta_{1})(\frac{1}{2}\mbox{ad}^{2}_{x}(\delta_{2}))^{l-1}
= \delta_{0}^{l},\end{align} which means that $ v_{2l-1},  v_{2l}
\leq l$. But on the other
 hand $ v_{2l-1},  v_{2l} \geq 0$ and
$$n(n+1) = v_{1} +...+  v_{2n} = \sum_{l=1}^{n}(v_{2l-1} + v_{2l}) \leq
\sum_{l=1}^{n}(l+ l) = n(n+1),$$ which implies the desired
equalities. $\Box $

\begin{definition} We will identify graded algebras
$\rm{Sym}(\rm{sl}_{2})$ and $k[z_{0}, z_{1}, z_{2}]$ using
$$z_{0}=e,\ \ z_{1}=-\frac{1}{2}h,\ \ z_{2}=-f,$$
where $deg\ z_{p} = p$. Then the canonical Kirillov-Kostant
Poisson structure takes the following form
\begin{align*}
\{ z_{0}, z_{1}\}  =  z_{0},  \ \ \{ z_{0}, z_{2}\}  =  2z_{1}, \
\ \{ z_{1}, z_{2}\}   = z_{2}.
\end{align*}

 Let ${\mathcal P} := k[z_{0},
z_{1}, z_{2}]/(z_{1}^{2} - z_{0}z_{2})$ be a graded polynomial
algebra of the Zariski closure of the unique not closed ${\rm
Ad}_{{\rm SL}_{2}}^{*}$-orbit in ${\rm sl}_{2}^{*}$ (which here
can be identified via the Killing form with the nilpotent cone in
${\rm sl}_{2}$).

The maximal  ideal of the vertex
$$\mathfrak{m} = (z_{0}, z_{1}, z_{2})$$
is a homogeneous Poisson ideal. We denote ${\mathcal
P}_{n}:={\mathcal P}/\mathfrak{m}^{n+1}$.
\end{definition}
\begin{theorem} There exists an isomorphism of graded polynomial Poisson
algebras
$${\mathcal P}_{n} \rightarrow {\mathcal P}({\mathcal O}_{n}),$$
$$(z_{0}, z_{1}, z_{2})\mapsto (y_{0}, y_{1}, y_{2}).$$
\end{theorem}

\textit{Proof.} By Lemma 4 we know that  the algebra ${\mathcal
P}({\mathcal O}_{n})$ is generated by $y_{p}$'s for $p = 0, ...,
2n$. By (9), (10) from the proof of Lemma 7 we have
$$\frac{1}{(2l)!}\mbox{ad}^{2l}_{x}(\delta_{2}^{l}) =  \delta_{0}^{l} =
\frac{1}{(2l)!}\mbox{ad}^{2l}_{x}(\delta_{2l}) ,$$
 $$\frac{1}{(2l-1)!}\mbox{ad}^{2l-1}_{x}(\delta_{1}\delta_{2}^{l-1}) = \delta_{0}^{l} =
\frac{1}{(2l-1)!}\mbox{ad}^{2l-1}_{x}(\delta_{2l-1}),$$ which
means that in ${\mathcal P}({\mathcal O}_{n})$
\begin{align*}
y_{2l}  =  y_{2}^{l}, \ \ y_{2l-1}  =  y_{1}y_{2}^{l-1}.
\end{align*}
Therefore  ${\mathcal P}({\mathcal O}_{n})$ is generated by $(
y_{0}, y_{1}, y_{2})$. By the multiplicative law (1) from Lemma 3
we have
$$\frac{1}{2}\mbox{ad}^{2}_{x}(\delta_{1}^{2}) = (\mbox{ad}_{x}(\delta_{1}))^{2} =
\delta_{0}^{2},$$
$$\frac{1}{2}\mbox{ad}^{2}_{x}(\delta_{0}\delta_{2}) = \delta_{0}
\frac{1}{2}\mbox{ad}_{x}^{2}(\delta_{2})=  \delta_{0}^{2}.$$
Subtracting these equations we get
$$\frac{1}{2}\mbox{ad}^{2}_{x}(\delta_{1}^{2} - \delta_{0}\delta_{2}) = 0,$$
which implies that in ${\mathcal P}({\mathcal O}_{n})$
\begin{align}
y_{1}^{2} - y_{0}y_{2} = 0. \end{align} This relation is also an
immediate consequence of the distinguished homomorphism.

Now using Lemma 4 we can write down the basis of  ${\mathcal
P}^{p}({\mathcal O}_{n})$ for every $p = 0, ..., 2n$

\begin{align}y_{2}^{l}, y_{0}y_{2}^{l}, ...,
y_{0}^{n-l}y_{2}^{l},\ \ \ \ \ {\rm for}\ \  p & =2l, \\
y_{1}y_{2}^{l}, y_{0}y_{1}y_{2}^{l}, ...,
y_{0}^{n-l-1}y_{1}y_{2}^{l},\ \ \ \ \ {\rm for}\ \  p & =2l+1.
\end{align}
Since $\delta_{0}^{n+1}=0$ we have
$$\frac{1}{p!}\mbox{ad}^{p}_{x}(\delta_{0}^{n+1-v_{p}}\delta_{p}) =
\delta_{0}^{n+1-v_{p}}\delta_{0}^{v_{p}} = \delta_{0}^{n+1} = 0,$$
which implies that
$$y_{0}^{n+1-v_{p}}y_{p} = 0.$$
By Lemma 7 this is equivalent to
\begin{align*}
y_{0}^{n+1-l}y_{2l-1}  =  0, \ \ y_{0}^{n+1-l}y_{2l}    =  0,
\end{align*}
and the latter by (11) is equivalent to
\begin{align*}
y_{0}^{n+1-l}y_{1}y_{2}^{l-1}  =  0, \ \ y_{0}^{n+1-l}y_{2}^{l} =
0.
\end{align*}
Again by (11) this implies that for all possible $p_{0}, p_{1},
p_{2}$ such that $p_{0} + p_{1} + p_{2} = n+1$
\begin{align}y_{0}^{p_{0}}y_{1}^{p_{1}}y_{2}^{p_{2}}=0.\end{align}
Looking at the basis (12)-(13) we see that  (11) and (14) form
together a complete system of relations on generators $ y_{0},
y_{1}, y_{2}$ which proves that the application of graded algebras
${\mathcal P}_{n}\rightarrow {\mathcal P}({\mathcal O}_{n})$,
$z_{p}\mapsto y_{p}$ for $p = 0, 1, 2$ is well defined and is an
isomorphism. The equality $\mbox{ad}_{x}(\delta_{1}) = \delta_{0}$
implies $\{ y_{0}, y_{1}\} = y_{0}$ and by the construction of
$\delta_{1} :=1/2\ \mbox{ad}_{x}(\delta_{2})$ we have $\{ y_{0},
y_{2}\} = 2y_{1}$. The following calculation
$$\mbox{ad}_{x}([\delta_{1}, \delta_{2}] - \delta_{2}) = [\mbox{ad}_{x}(\delta_{1}), \delta_{2}] +
[\delta_{1}, \mbox{ad}_{x}(\delta_{2})] -
\mbox{ad}_{x}(\delta_{2}) =
$$
$$ = [\delta_{0}, \delta_{2}] + [\delta_{1}, 2\delta_{1}] - [\delta_{0}, \delta_{2}] = 0$$
shows that $\{ y_{1}, y_{2}\} = y_{2}$. This relation is also an
immediate consequence of the distinguished choice of $\delta_{0},
\delta_{1}, \delta_{2}$. Therefore the above isomorphism of graded
algebras preserves the Poisson structure. $\Box $
\begin{corollary} ${\rm Spec}({\mathcal P}({\mathcal O}_{n}))$ admits a structure of a conical Poisson
scheme with a hamiltonian infinitesimal algebraic action of
$SL_{2}$, whose moment map is a closed embedding onto the $n$-th
infinitesimal neighborhood of the vertex in the nilpotent cone in
${\rm sl}_{2}^{*}$.
\end{corollary}

\begin{theorem} The distinguished homomorphism of filtered algebras
$${\mathcal U}({\rm sl}_{2})/(C - n(n + 2))\rightarrow {\mathcal D}({\mathcal O}_{n}),$$
is surjective with kernel $(e^{n+1})$.
\end{theorem}

{\em Proof.} Since the principal symbols $y_{0}, y_{1}, y_{2}$ of
the distinguished operators $\delta_{0}, \delta_{1}, \delta_{2}$
generate the associated graded algebra ${\mathcal P}({\mathcal
O}_{n})$ of the filtered algebra ${\mathcal D}({\mathcal O}_{n})$
then the distinguished operators generate ${\mathcal D}({\mathcal
O}_{n})$ as well, which proves the surjectivity of the
distinguished homomorphism. Since all maximal ideals in ${\mathcal
U}({\rm sl}_{2})$ of height two are of the form $(e^{n+1}, C - n(n
+ 2))$, $n\geq 0$ (\cite{cat}, Theorem 4.5 (ii)) this proves our
theorem. $\Box $

\vspace{3mm}
\paragraph{\bf 4. Principal symbols and inverse limits.} Note that in general the graded
Poisson algebra spanned by principal symbols is not functorial
with respect to a given commutative algebra. However graded
Poisson algebras spanned by principal symbols of the surjective
inverse system ${\mathcal O}_{n+1}\rightarrow {\mathcal O}_{n}$
(when applied object-wise) form a surjective inverse system
${\mathcal P}({\mathcal O}_{n+1})\rightarrow {\mathcal
P}({\mathcal O}_{n})$ as well, and we have
\begin{align}
\lim_{n}{\mathcal P}({\mathcal O}_{n})=k[[y_{0}, y_{1},
y_{2}]]/(y_{1}^{2}-y_{0}y_{2}).
\end{align}
On the other hand, for
${\mathcal O}=k[x]$ we have ${\mathcal P}({\mathcal O})=k[x_{0},
x_{1}]$, where $x_{0}$ and $ x_{1}$ are principal symbols of
operators $x$ and $\frac{{\rm d}}{{\rm d}x}$, respectively. Fix a
maximal ideal $\mathfrak{m}=(x_{0}, x_{1})\subset {\mathcal
P}({\mathcal O})$ and denote ${\mathcal P}({\mathcal
O})_{n}:={\mathcal P}({\mathcal O})/\mathfrak{m}^{n}$. Since
\begin{align}
\{ x_{0}, x_{1}\}=-1
\end{align}
this ideal and all its powers are not Poisson, so ${\mathcal
P}({\mathcal O})_{n}$ are not Poisson algebras. However
\begin{align}
\{ \mathfrak{m}^{m}, \mathfrak{m}^{n}\}\subset
\mathfrak{m}^{m+n-1},
\end{align}
hence the inverse limit
\begin{align}
\lim_{n}{\mathcal P}({\mathcal O})_{n}=k[[x_{0}, x_{1}]]
\end{align}
inherits a graded Poisson structure from ${\mathcal P}({\mathcal
O})$.

It is an interesting fact that these two complete graded Poisson
$\lim_{n}{\mathcal O}_{n}$-algebras can be compared in a unique
way. More presisely, we have the following theorem.

\begin{theorem}There exists the unique grading preserving Poisson homomorphism
of graded Poisson algebras making the following diagram
commutative
$$\begin{array}{ccccc}
\lim_{n}{\mathcal P}({\mathcal O})_{n}     &           &  \stackrel{\exists!}{\longleftarrow}  &           & \lim_{n}{\mathcal P}({\mathcal O}_{n}) \\
                                           & \nwarrow  &                                       & \nearrow  &                                        \\
                                           &           & \lim_{n}{\mathcal O}_{n} &            &
     \end{array}$$
     \end{theorem}

\textit{Proof:} We will look for a solution of the above problem
using the explicit presentation

$$\begin{array}{ccccc}
k[[x_{0}, x_{1}]] &           &  \stackrel{\exists?}{\longleftarrow}  &           & k[[y_{0}, y_{1}, y_{2}]]/(y_{1}^{2}-y_{0}y_{2}) \\
                  & \nwarrow  &                                       & \nearrow  &                                                 \\
                  &           &                 k[[x]]                &           &
     \end{array},$$
where the structure of $k[[x]]$-algebra is given by $x\mapsto
y_{0}$ and $x\mapsto x_{0}$, respectively.

As a homomorphism of $k[[x]]$-algebras the map $k[[y_{0}, y_{1},
y_{2}]]/(y_{1}^{2}-y_{0}y_{2})\rightarrow k[[x_{0}, x_{1}]]$
should be of the form
$$y_{0}\mapsto x_{0}, \ \ y_{1}\mapsto f_{1}, \ \ y_{2}\mapsto f_{2},$$
for some $f_{1}, f_{2}\in k[[x_{0}, x_{1}]]$ such that
\begin{align}
f_{1}^{2}=x_{0}f_{2},
\end{align}
which implies that
\begin{align}
f_{1}=x_{0}g, \ \ f_{2}=x_{0}g^{2},
\end{align}
for some $g\in k[[x_{0}, x_{1}]]$.

This is grading preserving if and only if $g$ is of degree one,
hence of the form
\begin{align}
g=hx_{1},
\end{align}
for some $h$ of degree zero, which means that $h\in
k[[x_{0}]]\subset k[[x_{0}, x_{1}]]$. Consequently
\begin{align}
f_{1}(x_{0}, x_{1})=x_{0}h(x_{0})x_{1}, \ \ f_{2}(x_{0},
x_{1})=x_{0}h(x_{0})^{2}x_{1}^{2}.
\end{align}

This is a Poisson homomorphism if and only if
\begin{align}
\{ x_{0}, f_{1}\} =x_{0},\ \ \{ x_{0}, f_{2}\} =2f_{1},\ \ \{
f_{1}, f_{2}\} =f_{2},
\end{align}
which by (16) and (22) is equivalent to
\begin{align}
h=-1.
\end{align}
All this means that there exists the unique solution to our
problem and it is of the form
\begin{align}
y_{0}\mapsto x_{0}, \ \ y_{1}\mapsto -x_{0}x_{1}, \ \ y_{2}\mapsto
x_{0}x_{1}^{2}. \ \ \Box
\end{align}

It is a very interesting homomorphism. First of all, it is the
completion of a (birational) resolution of singularity in the
category of graded Poisson algebras. Geometrically, this
resolution describes the following diagram in the category of
conical Poisson varieties
$$\begin{array}{ccccc}
{\rm T}^{*}\mathbb{A}^{1}     &           &  \stackrel{\exists!}{\longrightarrow} &           & Y \\
                              & \searrow  &                                       & \swarrow  &   \\
                              &           &            \mathbb{A}^{1}             &           &
     \end{array},$$
where $Y={\rm Spec}(k[y_{0}, y_{1},
y_{2}]/(y_{1}^{2}-y_{0}y_{2}))$ is a quadratic cone. Left hand
side morphism is a family of lines which is mapped into a family
of conics degenerating to a double line $(y_{0}=0, y_{1}^{2}=0)$
on the right hand side. This resolution contracts the line
$(x_{0}=0)$ to the vertex of the cone $Y$ and the scheme
theoretical pre-image of the double line is the line $(x_{0}=0)$
again.

\vspace{3mm} \paragraph{\textbf{5. Infinitesimal-birational
duality.}} Another remarkable property of this morphism deserves a
separate theorem. The following theorem means that the above
canonical resolution of singularity together with the Fourier
transformation establish a \textbf{duality} between the above
$\rm{SL}_{2}$-symmetries of principal symbols of the system of
infinitesimal neighborhoods of finite order of a point in the line
$\mathbb{A}^{1}$ (\textbf{infinitesimal} picture)  and
automorphisms of dual principal symbols induced by birational
automorphisms of the dual line $\hat{\mathbb{A}}^{1}$
(\textbf{birational } picture).

\begin{theorem}{\rm \textbf{``Magic Diamond Theorem"}}. The
following diagram  of canonical morphisms in the category of
conical Poisson varieties
$$\begin{array}{ccc}
                                              &  {\rm T}^{*}\mathbb{A}^{1}  &              \\
\stackrel{\stackrel{\mbox{\scriptsize{\rm resol.\ sing.}}}{}}{}\swarrow  &                             & \searrow\stackrel{\stackrel{\mbox{\scriptsize{\rm cl. Fourier\ tr.}}}{}}{}     \\
   Y          &                             &   {\rm T}^{*}\hat{\mathbb{A}}^{1}        \\
                                     \mbox{\scriptsize{\rm inf. moment\ m.}}\searrow &                             & \swarrow\mbox{\scriptsize{\rm bir. moment\ m.}}     \\
                                              &   {\rm sl}_{2}^{*}          &
                                   \end{array}$$
is commutative.
\end{theorem}

\textit{Proof:} The left hand side of this diagram  is already described in Corollary 5 and
Theorem 3. The composition on the left hand side reads as
$$k[z_{0}, z_{1},
z_{2}]\rightarrow k[x_{0}, x_{1}],$$
\begin{align}
z_{0}\mapsto x_{0},\ \ z_{1}\mapsto -x_{0}x_{1},\ \ z_{2}\mapsto
x_{0}x_{1}^{2}.
\end{align}
The classical Fourier transformation on the right hand side  reads as
$${\rm T}^{*}\hat{\mathbb{A}}^{1}\leftarrow {\rm T}^{*}\mathbb{A}^{1},$$
$$k[\hat{x}_{0}, \hat{x}_{1}] \rightarrow k[x_{0}, x_{1}],$$
\begin{align}
\hat{x}_{0}  \mapsto x_{1},\ \ \hat{x}_{1} \mapsto -x_{0}.
\end{align}

The infinitesimal form of the  rational action of the algebraic
group $\rm{SL}_{2}$ on ${\rm T}^{*}\hat{\mathbb{A}}^{1}$
canonically induced from the canonical action of $\rm{SL}_{2}$ on
$\hat{\mathbb{A}}^{1}$ by birational automorphisms
$$\left( \begin{array}{cc}
         a & b \\
         c & d
         \end{array}\right)\cdot (\hat{x}_{0}, \hat{x}_{1})=
         \left( \frac{a\hat{x}_{0}+b}{c\hat{x}_{0}+d},(c\hat{x}_{0}+d)^{2}\hat{x}_{1}\right)$$
reads as follows
\begin{align}
e\mapsto &- \frac{\partial}{\partial \hat{x}_{0}} = \{-\hat{x}_{1}, -\},\\
h\mapsto & -2\hat{x}_{0}\frac{\partial}{\partial
\hat{x}_{0}}+2\hat{x}_{1}\frac{\partial}{\partial \hat{x}_{1}}=\{-2\hat{x}_{0}\hat{x}_{1}, -\},\\
f\mapsto & \hat{x}_{0}^{2}\frac{\partial}{\partial
\hat{x}_{0}}-2\hat{x}_{0}\hat{x}_{1}\frac{\partial}{\partial
\hat{x}_{1}}=\{\hat{x}_{0}^{2}\hat{x}_{1}, -\}.
\end{align}
Therefore the moment map associated with this action on ${\rm
T}^{*}\hat{\mathbb{A}}^{1}$ is of the form
$$k[z_{0}, z_{1},
z_{2}]\rightarrow k[\hat{x}_{0}, \hat{x}_{1}],$$
$$z_{0}\mapsto -\hat{x}_{1},\ \ z_{1}\mapsto -\hat{x}_{0}\hat{x}_{1},\ \ z_{2}\mapsto -\hat{x}_{0}^{2}\hat{x}_{1}.$$
It is easy to see that its composition with the classical Fourier
transformation (27) gives the same  as (26), what proves the
commutativity of the above square. $\Box $

\vspace{3mm}
\paragraph{\bf Remark.}  The moment map on ${\rm
T}^{*}\hat{\mathbb{A}}^{1}$ is the restriction of the canonical
moment map on ${\rm T}^{*}\hat{\mathbb{P}}^{1}$ associated with
the action of the automorphism group of $\hat{\mathbb{P}}^{1}$.
One has the following commutative diagram
$$\begin{array}{ccccc}
                        & {\rm T}^{*}\hat{\mathbb{P}}^{1} & \stackrel{\cong}{\longleftarrow}      & \tilde{\mathcal{N}} & \\
  \mbox{\tiny{moment m.}} & \downarrow                      &                                   &\downarrow           & \mbox{\tiny{Springer res.}} \\
                        & \rm{sl}_{2}^{*}                 & \stackrel{\rm{Killing\ f.}}{\longleftarrow} & \mathcal{N} &
  \end{array}$$
where the right-hand vertical morphism is the Springer resolution
of the nilpotent cone $\mathcal{N}\subset \rm{sl}_{2}$ \cite{jan}.
Therefore, by Theorem 4, our canonical resolution as in Theorem 3
is, up to the classical Fourier transformation and the isomorphism provided
by the Killing form, essentially nothing but the Springer
resolution.

\vspace{3mm} \paragraph{\textbf{6. Quantization of the infinitesimal-birational
duality.}} The following theorem shows and explains the quantum nature of the above duality.

\begin{theorem}{\rm \textbf{``Quantized Magic Diamond Theorem"}}. There exists the
following commutative diagram  in the category of
almost commutative algebras
$$\begin{array}{ccc}
                                              &  \mathcal{D}(\mathbb{A}^{1})  &              \\
\stackrel{\stackrel{\mbox{\scriptsize{\rm q. resol.\ sing.}}}{}}{}\nearrow  &                             & \nwarrow\stackrel{\stackrel{\mbox{\scriptsize{\rm q. Fourier\ tr.}}}{}}{}     \\
  \mathcal{U}({\rm sl}_{2})/(C-n(n+2))   &                             &    \mathcal{D}(\hat{\mathbb{A}}^{1})        \\
                                     \mbox{\scriptsize{\rm q. inf. moment\ m.}}\nwarrow &                             & \nearrow\mbox{\scriptsize{\rm q. bir. moment\ m.}}     \\
                                              &   \mathcal{U}({\rm sl}_{2})          &
                                   \end{array}$$
which is a quantization of the Magic Diamond.
\end{theorem}
\textit{Proof:} The only nonobvious arrows in this diagram are: the quantized birational moment map on the right hand side and the quantized resolution of singularities on the left hand side.  The first one is canonically obtained from the birational action of ${\rm SL}_{2}$ on  $\hat{\mathbb{A}}^{1}$ as above, i.e.
\begin{align}
e\mapsto -\frac{d}{d\hat{x}},\ \ \ h\mapsto -2\hat{x}\frac{d}{d\hat{x}},\ \ \  f\mapsto \hat{x}^{2}\frac{d}{d\hat{x}}.
\end{align}
Composing with the quantum Fourier transformation we get
\begin{align}
e\mapsto x,\ \ \ h\mapsto 2x\frac{d}{dx}+2,\ \ \  f\mapsto -x^{2}\frac{d^{2}}{dx^{2}}-2\frac{d}{dx}.
\end{align}
Note that these operators look exactly like the operators (8) inducing the distinguished differential operators on the $n$-th infinitesimal neighborhood of a point in $\mathbb{A}^{1}$, except the strange value $n=-2$. Nevertheless, they satisfy relation $h^{2}+2(ef+fe)=n(n+2)$ for $n=-2$, which means that the composite homomorphism factorizes as in the left hand side of the diagram. This proves the commutativity of this diagram. Passing to the diagram of associated graded Poisson algebras is equivalent to  replacing all differential operators (32) by their principal symbols and all differential operators (31) by their classes modulo differential operators with polynomial coefficients of lower degree. It is easy to see that we reobtain  the Magic Diamond.   $\Box$

\vspace{3mm}
\paragraph{\bf Remark.}  The strange value $n=-2$ in the above proof justifies our definition of  quantization as Morita equivalence class. Namely, for $n\geq 0$ the primitive quotient $\mathcal{U}({\rm sl}_{2})/(C-n(n+2))$ is a  representative of a quantization of the nilpotent cone admitting closed embedding (surjective homomorphism of algebras) of  a representative of  a quantization of $n$-th infinitesimal neighborhood of the vertex (see Corollary 5 and Theorem 2). But by the Beilinson-Bernstein Theorem these primitive quotients are Morita equivalent for all $n\in\mathbb{Z}$. This means that  this quantization can  be computed as the class of any primitive quotient $\mathcal{U}({\rm sl}_{2})/(C-n(n+2))$, even for $n=-2$ which has no any immediate geometric meaning.

\end{document}